\documentclass[12pt]{amsart}
\usepackage{amssymb}
\usepackage{amsthm}
\usepackage{graphicx} 
\usepackage{textcomp}
\usepackage{extpfeil}

\newtheorem*{thm}{Theorem}
\newtheorem*{claim_V_1}{Claim \ref{claim_V_1}}

\newtheorem*{E1s1}{Theorem \ref{E1s1}}
\newtheorem{theorem}{Theorem}[section]

\newtheorem{lemma}[theorem]{Lemma}
\newtheorem{corollary}[theorem]{Corollary}
\newtheorem{question}{Question}

\newtheorem{claim}{Claim}
\newtheorem*{conjecture}{Conjecture}
\title[RANK GRADIENT IN CO-FINAL TOWERS]{RANK GRADIENT IN CO-FINAL TOWERS OF CERTAIN KLEINIAN GROUPS}
\author{Darlan Gir\~ao}

\begin{document}
\maketitle
\footnote{Author supported by CAPES/Fulbright Grant BEX 2411/05-9}
\begin{abstract}
We prove that if the fundamental group of an orientable finite volume  hyperbolic 3-manifold has finite index in the reflection group of a right-angled ideal polyhedra in $\mathbb{H}^3$ then it has a co-final tower  of finite sheeted covers with positive rank gradient. The manifolds we provide are also known to have co-final towers of covers with zero rank gradient.  
\end{abstract}

\section{INTRODUCTION}

Let $G$ be a finitely generated group. The \textit{rank of $G$} is the minimal cardinality  of a generating set, and is denoted by  $\text{rk}(G)$.  If $G_j$ is a finite index subgroup of $G$, the Reidemeister-Schreier process (\cite{LS}) gives an upper bound on the rank of $G_j$.	 
$$\text{rk}(G_j)-1\leq [G:G_j](\text{rk}(G)-1)$$  
Recently Lackenby introduced the notion of \textit{rank gradient} (\cite{La1}).
 Given a finitely generated group $G$ and a collection $\{G_j\}$ of finite index subgroups,  the \textit{rank gradient} of  the pair $(G,\{G_j\})$ is defined by 
$$\text{rgr}(G,\{G_j\})=\lim_{j\rightarrow\infty}\frac{\text{rk}(G_j)-1}{[G:G_j]}$$   
We say that the  collection of finite index subgroups $\{G_j\}$ is \textit{co-final} if $\cap_j G_j=\{1\}$, and we call it a  \textit{tower} if $G_{j+1}<G_j$.

In some particular cases it is easy to determine rank gradient, for example:  
\begin{itemize}
\item[(1)] When $G$ is a free group, the rank gradient of any pair $(G,\{G_j\})$ is positive. 
\item[(2)] The same is true if $G$ is the fundamental group of a closed surface $S$ with $\chi(S)<0$; 
\item[(3)] If $G\xtwoheadrightarrow{}F_2$, where  $F_2$ is the free group on two generators then, using (1), one can find a tower (not co-final)  of subgroups with positive rank gradient;  
\item[(4)] If $G$ is virtually abelian or if $G$ is the fundamental group of a virtually fibered 3-manifold then there are towers  with zero rank gradient. In the latter case we consider the subgroups coming from the cyclic covers of the fibered manifold.    
\item[(5)] $\text{SL}(n,\mathbb{Z}), n>2$, has zero rank gradient with respect to towers of congruence subgroups (\cite{Ti}, \cite{La1}). 
\end{itemize}
However, determining the rank gradient of a co-final tower is very hard in general.  For example, the following question is the motivation for this note:

\begin{question} Does there exist  a torsion free finite covolume Kleinian group $G$ with a co-final tower $\{G_j\}$  such that  $\operatorname{rgr}(G,\{G_j\})>0$. 
\end{question}
\noindent
The main result of this note provides infinitely many such examples. To state it we introduce some notation.

If $M_1$ is an orientable finite volume hyperbolic 3-manifold, we call the family of covers $\{M_j\longrightarrow M_1\}$ \textit{co-final} (resp. a \textit{tower}) if $\{\pi_1(M_j)\}$ is co-final (resp. a tower). By rank gradient of the the pair $(M_1,\{M_j\})$,  $\text{rgr}(M_1,\{M_j\})$, we mean the rank gradient of $(\pi_1(M_1),\{\pi_1(M_j)\})$. 

\begin{E1s1}
Let $M_1$ be an orientable finite volume hyperbolic 3-manifold whose fundamental group has finite index in the reflection group of a totally geodesic right-angled ideal polyhedron $P_1$ in $\mathbb{H}^3$. Then there exists a  co-final tower of finite sheeted covers $\{M_j\longrightarrow M\}$ with positive rank gradient.
\end{E1s1}
  
This theorem relates to the work of  Ab\'ert and Nikolov (\cite{AN}), and in particular to a question about \textit{cost of group actions} (\cite{Ga}). 

\begin{question}\label{cost}
Let $G$ be finitely generated and $\{G_j\}$ be a co-final tower  of  normal subgroups of $G$.  
Does $\operatorname{rgr}(G,\{G_j\})$ depend on the tower  $\{G_j\}$? 
\end{question}

Our result provides  negative evidence for this question. If one could improve Theorem \ref{E1s1} by finding a co-final tower $\{M_j\longrightarrow M_1\}$ of regular covers  with positive rank gradient, then we claim it would also be possible to find one with zero rank gradient.  In fact,  Agol proved in \cite{Ag} that if the fundamental group of a 3-dimensional manifold satisfies an algebraic condition, called RFRS,  then it virtually fibers. He also proved in \cite{Ag} that the manifolds of the type considered in Theorem \ref{E1s1} are virtually RFRS. Therefore, given $M_1$ as in Theorem \ref{E1s1}, it is possible to find a tower  $\{\Gamma_j\}$ with $\text{rgr}(\pi_1(M_1),\{\Gamma_j\})=0$. By taking the \textit{core} of $\Gamma_j$ in $\pi_1(M_1)$ (i.e., $\text{core}(\Gamma_j)=\displaystyle\cap_{g\in\pi_1(M_1)}g\Gamma_jg^{-1}$), one sees that the tower of normal subgroups $\{\text{core}(\Gamma_j)\}$ has zero rank gradient. The desired co-final tower with zero rank gradient would be given by   $\{\pi_1(M_j)\cap\text{core}(\Gamma_j)\}$. 

The main idea of the proof of Theorem \ref{E1s1} is as follows: given $P_1$ as in the theorem, construct a collection of polyhedra $\{P_j\}$ whose reflection groups have finite  index $2^{j-1}$ in the reflection group of $P_1$. If one is given an orientable hyperbolic 3-manifold $M_1$  whose fundamental group has finite index in the reflection group of  $P_1$  then $M_1$ has at least as many cusps as the number of vertices of $P_1$.  We may find manifold covers $M_j\longrightarrow M_1$ so that $M_j$ is  a $2^{j-1}$-sheeted covering and has at least as many cusps as the number of ideal vertices of $P_j$. We then show that the $P_j$ can be chosen so that the number of its vertices is of the same magnitude as $2^{j}$. 

The paper will be organized as follows: section 2  sets up notation and we recall a characterization of right-angled ideal polyhedra using Andreev's theorem (\cite{An}). We then show how the construction of the family $\{P_j\}$ will be done. In section 3 we prove  Theorem \ref{E1s1}. Section 4 contains all the technical results we need to estimate $\text{rk}(\pi_1(M_j))$. In section 5 we show how to construct $\{P_j\}$ so that the family $\{M_j\}$ is co-final. The idea for this appears in \cite{Ag} (Theorem 2.2) and we include  a proof  here for completeness. Section 6 contains some final remarks and further questions.

\section*{ACKOWLEDGEMENTS}
I am  very grateful to my thesis advisor, Alan Reid, for his extraordinary guidance and unwavering support. I am also thankful  to Ian Agol, Grant Lakeland and Mark Norfleet for helpful  conversations. The author was partially supported by CAPES/Fulbright Grant BEX 2411/05-9.    

\section{SET UP} 
An \textit{abstract polyhedron} $\mathcal{P}_1$ is a cell complex on $S^2$ which can be realized by a convex
Euclidean polyhedron.   A \textit{labeling} of  $\mathcal{P}_1$ is a map 
$$ \Theta: \text{Edges}(\mathcal{P}_1)\longrightarrow  (0,\pi/2]$$
The pair $(\mathcal{P}_1,\Theta)$ is a labeled abstract
polyhedron. A labeled abstract polyhedron is said to be \textit{realizable} as a hyperbolic
polyhedron if there exists a hyperbolic polyhedron $P_1$ such that there is a label preserving
graph isomorphism between the 1-skeleton of $P_1$ with edges labeled by dihedral angles
and the 1-skeleton of $\mathcal{P}_1$ with edges labeled by $\Theta$.

Let $P_1$ be a totally geodesic right-angled ideal polyhedron in $\mathbb{H}^3$ (that is, faces of $P_1$ are contained in hyperplanes and all vertices of $P_1$ lie in the boundary at infinity $S^2_{\infty}$, where we here we consider the ball model for $\mathbb{H}^3$).
We consider  the 1-skeleton of $P_1$ as a graph $\Gamma_1 \subset S^2$ with labels $\theta_e=\pi/2$. Let $\Gamma_1^*$ be its dual graph.  A \textit{$k$-circuit} is a simple closed curve composed of $k$ edges in $\Gamma_1^*$. A \textit{prismatic $k$-circuit} is a $k$-circuit $\gamma$ 
 so that no two edges of $\Gamma_1$ which correspond to edges traversed by $\gamma$ share a vertex.  Andreev's   theorem for 
right-angled ideal  polyhedra  in $\mathbb{H}^3$ (\cite{An}, see also \cite{At})  can be stated as: 
\begin{theorem} Let $\mathcal{P}_1$ be an abstract polyhedron. Then $\mathcal{P}_1$ is realizable as a right-angled ideal polyhedron $P_1$ if and only if 
 \begin{itemize}
 \item [(1)] $P_1$ has at least 6 faces;
 \item [(2)] Vertices have valence 4;
 \item [(3)] For any triple of faces of $P_1$, $(f_i, f_j , f_k)$, such that $f_i\cap f_j$ and $f_j\cap f_k$ are edges of $P_1$ with distinct endpoints, $f_i\cap f_k=\emptyset$;.
 \item [(4)] There are no prismatic 4-circuits.
 \end{itemize}
\end{theorem}

The above theorem implies that the 1-skeleton of $P_1$ is a 4-valent graph. The faces can therefore be checkerboard colored. 
Reflecting $P_1$ along a face $f_1$ gives a polyhedron $P_2$ which is
also right-angled, ideal and totally geodesic with checkerboard colored faces (see figure below).  We construct a sequence of polyhedra $P_1, P_2,...,P_j,...$  recursively, whereby
$P_{j+1}$ is obtained from $P_j$ by reflection along a face $f_j$. The faces of $P_{j+1}$ are colored accordingly with the coloring of the faces of $P_{j}$. 

The notation for the remainder of the paper is as follows: the number of vertices in the face $f_j$ is denoted by $S_{f_j}$ and $\phi_{f_j}$ denotes the reflection along $f_j$. $B_j$ and $W_j$  represent the maximal number of ideal vertices on a black or white face of the polyhedron $P_j$, respectively.  $V_j$ denotes the total number of vertices on $P_j$.  

Throughout,  the construction of the polyhedra $P_j$ will be done in an
alternating fashion with respect to the color of the faces: $P_{2j}$ is obtained
from $P_{2j-1}$ by refection along a black face and $P_{2j+1}$ is obtained from
$P_{2j}$ by reflection along a white face. 
 
\begin{figure}[h]\label{oct}
\includegraphics[scale=.3]{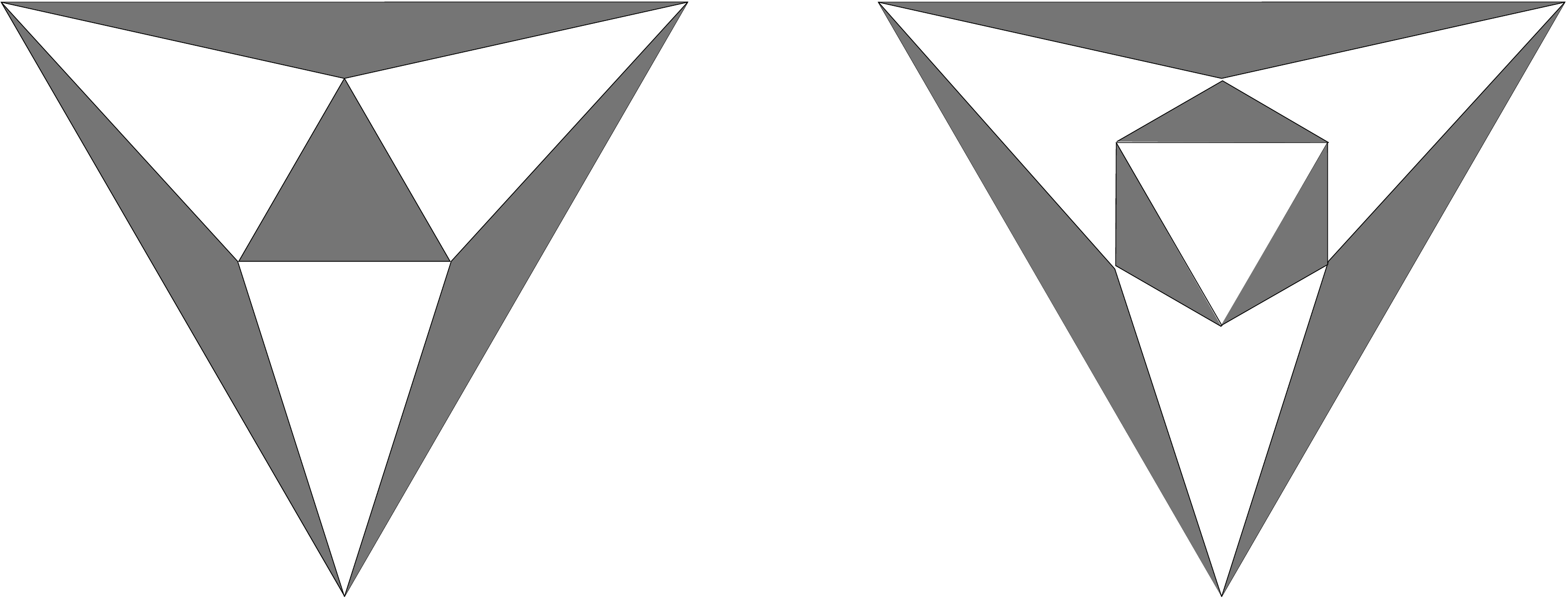}
\caption{Polyhedron $P_1$ reflected along central black face yields $P_2$}
\end{figure}

\section{ MAIN THEOREM}
In this section we prove:
\begin{theorem}\label{E1s1}
Let $M_1$ be an orientable finite volume hyperbolic 3-manifold whose fundamental group has finite index in the reflection group of a  right-angled ideal polyhedron $P_1$ in $\mathbb{H}^3$. Then there exists a co-final tower of finite sheeted covers $\{M_j\longrightarrow M\}$ with positive rank gradient.
\end{theorem}

Our construction of the family $\{M_j\}$ was inspired by  the proof of Theorem 2.2 of Agol's paper (\cite{Ag}). The proof that this family can be made co-final is given in  section 5 (following \cite{Ag}).   

\begin{proof}[Proof of Theorem \ref{E1s1}]
Consider the family of polyhedra $\{P_j\}$ obtained from $P_1$ as decribed above. Denote by $G_j$ the reflection group of $P_j$ and observe that $G_{j+1}$ is a subgroup of $G_j$ of index 2.   $G_1$  acts on $\mathbb{H}^3$ with fundamental domain $P_1$. The orbifold $\mathbb{H}^3/G_1$ is non-orientable, and may be viewed as  $P_1$ with its faces mirrored. The  singular locus is the 2-skeleton of $P_1$. Each  ideal vertex of $P_1$ corresponds to a cusp of $\mathbb{H}^3/G_1$. 

Let  $M_1$ be an orientable  cusped hyperbolic 3-manifold such that $\pi_1(M_1)$  has finite index in $G_1$. Let $M_j\longrightarrow M_1$ be the cover of $M_1$ whose fundamental group is  $\pi_1(M_j)=\pi_1(M_1)\cap G_j$. Since $[G_j:G_{j+1}]=2$, we must have $[\pi_1(M_j):\pi_1(M_{j+1})]\leq 2$. Also note that since $\text{vol}(P_j)=2^{j-1}\text{vol}(P_1)$, for all but finitely many $j$ (at most $[G_1:\pi_1(M_1)]$) we must have $[\pi_1(M_j):\pi_1(M_{j+1})]=2$. We may thus assume that $[\pi_1(M_j):\pi_1(M_{j+1})]=2$.  By mirroring the faces of $P_j$, it may be regarded as  a non-orientable finite volume orbifold (as described before). This implies that  $M_j\longrightarrow P_j$ is an orientable finite sheeted cover for $j=1,2,...$.

 Note that $[\pi_1(M_1):\pi_1(M_j)]=2^{j-1}$. Thus to show that the family $\{M_j\longrightarrow M_1\}$ has positive rank gradient we will establish that  $\text{rk}(\pi_1(M_j))$ grows with the same magnitude as $2^j$. 

By \textquotedblleft half lives half dies\textquotedblright, an easy lower bound on the rank of the fundamental group of an orientable finite volume  hyperbolic 3-manifold is the number of its cusps. Since the cusps of $P_j$ correspond to its ideal vertices and the number of cusps does not go down under finite sheeted covers, it must be that $M_j$ has at least as many cusps as the number of ideal vertices of $P_j$. 

Recall that  $B_j$ and $W_j$ are the maximal number of ideal vertices on a black or white face of the polyhedron $P_j$, respectively, and  $V_j$ is the total number of vertices on $P_j$. The claims below (proved in section 4) gives us the estimates we need for $V_j$ in terms of $V_1$, $B_1$ and $W_1$. 
  
\begin{claim}\label{claim_V_1}
$V_1\geq B_1+W_1-1$
\end{claim}

\begin{claim}\label{claim_V_j}
For any $j\geq 6$, 
$$V_j\geq 2^{j-1}V_1-2^{j-1}(B_1+W_1)+2^{j-1}+2^{j-2}$$
\end{claim}

Given these, we argue as follows:
$$\text{rgr}(M_1,\{M_j\})=\lim_{j\rightarrow\infty}\dfrac{\text{rk}(\pi_1(M_j))-1}{[\pi_1(M_1):\pi_1(M_j)]}\geq$$ $$\lim_{j\rightarrow\infty}\dfrac{V_j-1}{2^{j-1}}\geq\lim_{j\rightarrow\infty}\dfrac{
2^{j-1}V_1-2^{j-1}(B_1+W_1)+2^{j-1}+2^{j-2}-1}{2^{j-1}}\geq$$

$$\lim_{j\rightarrow\infty}\dfrac{2^{j-1}(B_1+W_1-1)-2^{j-1}(B_1+W_1)+2^{j-1}+2^{j-2}-1}{2^{j-1}}\geq$$ $$\lim_{j\rightarrow\infty}\dfrac{2^{j-2}-1}{2^{j-1}}=\dfrac{1}{2}$$ which proves the theorem.
\end{proof}

\section{LOWER BOUNDS ON NUMBER OF IDEAL VERTICES OF $P_j$}
  
We now proceed to prove Claims \ref{claim_V_1} and \ref{claim_V_j}. This requires several preliminary results.

\begin{lemma}\label{eqclass} Let $P_{j+1}$ be obtained from $P_j$ by
reflection along a face $f_j$.  Then $V_{j+1}= 2V_j-S_{f_j}$. 
\end{lemma}

\begin{proof}
Here we abuse notation and write $v\in f_j$ if $v$ is an ideal vertex of the face $f_j$ and write $v\notin f_j$ otherwise. 
Note that if  $v\notin f_j$, then $v$ yields two vertices on $P_{j+1}$, namely, $v$ and $\phi_{f_j}(v)$. If $v\in f_j$, then it yields a single vertex ($v$ itself).

If $v\notin f_j$, then, by the observation above, $v$ yields two
ideal vertices on $P_{j+1}$.  Since  a total of $S_{f_j}$  ideal vertices lie in $f_j$ 
and  $V_j-S_{f_j}$ do not, it must be that 
that  
$$V_{j+1}=2(V_j-S_{f_j})+S_{f_j}=2V_j-S_{f_j}$$
\end{proof}

Recall also that the construction of the family of polyhedra $\{P_j\}$ is made  in
an alternating fashion with respect to the color of the faces: $P_{2j}$ is
obtained from $P_{2j-1}$ by refection along a black face and $P_{2j+1}$ is
obtained from $P_{2j}$ by reflection along a white face.

\begin{corollary}\label{V_est} For $j\geq 1$
\begin{itemize}
\item[(1)] $V_{2j}\geq 2V_{2j-1}-B_{2j-1}$ 
\item [(2)] $V_{2j+1}\geq 2V_{2j}-W_{2j}$
\end{itemize}

\end{corollary}

\begin{proof} $P_{2j}$ is obtained from $P_{2j-1}$ by refection along a black face $f_{2j-1}$, thus $S_{f_{2j-1}}\leq B_{2j-1}$. By the lemma, $V_{2j}=2V_{2j-1}-S_{f_{2j-1}}$ and therefore $V_{2j}\geq 2V_{2j-1}-B_{2j-1}$ . The second inequality is similar.
\end{proof}

 With  the notation  established above we now find  lower bounds for the
$V_j$ in terms of $V_1, B_1$ and $W_1$. First we need to find upper
bounds for $B_j$ and $W_j$ in terms of $B_1$ and $W_1$. To do this in a way that will fit our purposes we establish two properties of the family $\{P_j\}$. As before, denote by $\phi_{f_j}$ the reflection along the face $f_j$. 

\begin{lemma}\label{s_j est}
\begin{itemize}
\item[(1)] If $P_j$ is reflected  along a white (resp. black) face $f_j$, all
black faces $f_*$ (resp. white faces $f_*$) adjacent to $f_j$ yield new black faces
$\tilde{f_*}$ (resp. white faces $\tilde{f_*}$) on $P_{j+1}$. The number
$S_{\tilde{f_*}}$ (resp. $S_{\tilde{f_*}}$) of ideal vertices on $\tilde{f_*}$
(resp. $\tilde{f_*}$) is  $2S_{f_*}-2$ (resp. $2S_{f_*}-2$).
\item[(2)]   A face $f_*$  not adjacent to $f_j$ yield two new faces, $f_*$ itself and $\phi_f(f_*)$, both with $S_{f_*}$ vertices.
\end{itemize}
\end{lemma}

\begin{proof}  For the first property, reflecting $f_*$  along $f_j$ gives a face $\phi_{f_j}(f_*)$ in $P_{j+1}$ adjacent to $f_*$. The dihedral angle between $f_*$ and $\phi_f(f_*)$
is $\pi$. Thus, on $P_{j+1}$, they correspond to  a single  face denoted by
$\tilde{f_*}$. The number of ideal vertices on $\tilde{f_*}$ is exactly
$2S_{f_*}-2$. The second property should be  clear.  See figure 1 for an ilustration of these properties.
\end{proof}
As an immediate consequence we have
\begin{corollary}\label{BW_est}
\begin{itemize}
\item[]
\item [(1)] 
$\begin{cases} 
B_{2j}=B_{2j-1}\\
W_{2j}\leq 2W_{2j-1}-2
\end{cases}$
\item [(2)] 
$\begin{cases}
B_{2j+1}\leq 2B_{2j}-2\\
W_{2j+1}=W_{2j}
\end{cases}$

\end{itemize}
\end{corollary}

We are now in position to estimate the values $B_j$ and $W_j$ in terms of $B_1$ and $W_1$.
\begin{theorem}
With the notation as before we have
\begin{itemize}
\item[(1)] $W_{2j+1}=W_{2j}\leq 2^jW_1-\displaystyle\sum_{l=1}^{j}2^l$
\item[(2)] $B_{2j+2}=B_{2j+1}\leq 2^jB_1-\displaystyle\sum_{l=1}^{j}2^l$
\end{itemize}
\end{theorem}

\begin{proof}  We procced by induction. By corollary \ref{BW_est} these statements are true for $j=1$.  Suppose it is also true for $j\leq n$. We now want to estimate $B_{2n+3}=B_{2n+4}$ and
$W_{2n+2}=W_{2n+3}$.  The hypothesis is that 
 $$W_{2j+1}=W_{2j}\leq 2^nW_1-\displaystyle\sum_{l=1}^{n}2^l$$
 $$B_{2n+2}=B_{2n+1}\leq 2^nB_1-\displaystyle\sum_{l=1}^{n}2^l$$  
$P_{2n+2}$ is obtained from $P_{2n+1}$
by reflection along a black face, denoted by $f$. 
White faces on $P_{2n+1}$ adjacent to $f$ yield new white
faces on $P_{2n+2}$ with at most $2W_{2n+1}-2$ vertices, by Corollary \ref{BW_est}. But 
$$2W_{2n+1}-2\leq 
2[2^nW_1-\displaystyle\sum_{l=1}^{n}2^l]-2=2^{(n+1)}W_1-\displaystyle\sum_{l=1}^
{n+1}2^l$$ which gives the desired result for $W_{2n+2}$ and $W_{2n+3}$.  Finally, $P_{2n+3}$
is obtained from $P_{2n+2}$ by a reflection
along a white face, again denoted by $f$.   Since black faces of
$P_{2n+2}$ have at most $B_{2n+2}(=B_{2n+1})$ vertices, black faces of $P_{2n+3}$ will have at most $2B_{2n+1}-2$ vertices, again by corollary \ref{BW_est}. But 
$$ 2B_{2n+1}-2\leq 
2[2^nB_1-\displaystyle\sum_{l=1}^{n}2^l]-2=2^{(n+1)}B_1-\displaystyle\sum_{l=1}^
{n+1}2^l$$ vertices. This establishes the result for $B_{2n+3}$ and $B_{2n+4}$.      
\end{proof}

\begin{theorem}\label{eqest}
With the notation as before, and for $j\geq 3$, 
\begin{itemize}
    \item[(1)] $V_{2j}\geq
2^{2j-1}V_1-B_1\displaystyle\sum_{l=j-1}^{2j-2}2^l - W_1\displaystyle\sum_{l=j}^{2j-2}2^l + \displaystyle\sum_
{l=j+2}^{2j-1}2^l+2^j+2$
     \item[(2)] $V_{2j+1}\geq 
2^{2j}V_1-B_1\displaystyle\sum_{l=j}^{2j-1}2^l - W_1\displaystyle\sum_{l=j}^{2j-1}2^l + \displaystyle\sum_{l=j+2}^{2j}
2^l+2$
 \end{itemize}
\end{theorem}

\begin{proof}
Lower bounds estimates for $V_1,..., V_7$  are found recursively.
$V_1$, $V_2$, $V_3$, $V_4$ and $V_5$ do not fit these formulas but $V_6$ and $V_7$  do. The statement is then true for $j=3$. We now proceed by induction, using   the previous proposition and corollary \ref{V_est}.   Suppose it is true for $j\leq n, n\geq 3$. We want to show this
implies true for $j=n+1$. By corollary \ref{V_est}, $V_{2n+2}\geq
2V_{2n+1}-B_{2n+1}$. The hypothesis is that    
$$V_{2n+1}\geq 
2^{2n}V_1-B_1\displaystyle\sum_{l=n}^{2n-1}2^l - W_1\displaystyle\sum_{l=n}^{2n-1}2^l + \displaystyle\sum_{l=n+2}^{2n}
2^l+2$$  
We also know that $$B_{2n+1}\leq 2^nB_1-\displaystyle\sum_{l=1}^{n}2^l$$
Thus 

$$V_{2n+2}\geq 2V_{2n+1}-B_{2n+1}\geq$$ 

$$2[2^{2n}V_1-B_1\displaystyle\sum_{l=n}^{2n-1}2^l - W_1\displaystyle\sum_{l=n}^{2n-1}2^l + \displaystyle\sum_{l=n+2}^{2n}
2^l+2]-[2^nB_1-\displaystyle\sum_{l=1}^{n}2^l]=$$

$$2^{2n+1}V_1-B_1\displaystyle\sum_{l=n}^{2n-1}2^{l+1} - W_1\displaystyle\sum_{l=n}^{2n-1}2^{l+1} + \displaystyle\sum_{l=n+2}^{2n}
2^{l+1}+2^2 + \displaystyle\sum_{l=1}^{n}2^l=$$

$$2^{2n+1}V_1-B_1\displaystyle\sum_{l=n}^{2n}2^{l} - W_1\displaystyle\sum_{l=n+1}^{2n}2^{l} + \displaystyle\sum_{l=n+3}^{2n+1}
2^{l}+2^{n+1}+2$$
\\
which establishes $(1)$  for $2(n+1)=2n+2$.

We use the exact same idea to and the estimate for $V_{2n+2}$ to establish (2) for
$2(n+1)+1=2n+3$.
\end{proof}

\begin{corollary}\label{V_j_est} 
For any $j\geq 6$, 
$$V_j\geq 2^{j-1}V_1-2^{j-1}(B_1+W_1)+2^{j-1}+2^{j-2}$$

\end{corollary}
\noindent
Hence Claim \ref{claim_V_j} in the proof of Theorem \ref{E1s1} is proved. We now prove
 
\begin{claim_V_1}
$V_1\geq B_1+W_1-1$
\end{claim_V_1}

\begin{proof}
Let $f_b$ and $f_w$ be black and white faces of $P_1$ with maximal number of vertices, i.e., $S_{f_b}=B_1$ and $S_{f_w}=W_1$. \\ 

 \textit{Case 1: The faces $f_b$ and $f_w$ are not adjacent} 

Here we get $V_1\geq B_1+W_1$ and the claim follows.\\  

\textit{Case 2: The faces $f_b$ and $f_w$ are adjacent}. 

Since $f_b$ and $f_w$ share exactly 2 vertices we see that $V_1\geq B_1+W_1-2$. Suppose we have equality. Then every vertex of  $P_1$ must be a vertex of either $f_b$ or $f_w$. Recall that we can visualize the 1-skeleton of $P_1$ as lying in $S^2$. Label the vertices of $P_1$ by $\{v_1,..., v_k\}$. The assumption is that all these vertices lie in the boundary of the disk $D=\overline{(f_b\cup f_w)}\subset S^2$. By  Andreev's theorem, $P_1$ has at least 6 faces, every face  is at least 3-sided and all vertices  are 4-valent. Denoting by $F_1$ and $E_1$ the number of faces and edges of $P_1$ respectively we have the relation $V_1-E_1+F_1=2$. Since  vertices are 4-valent we also have  $E_1=2V_1$. From these relations and $F_1\geq 6$, we get $V_1\geq 4$.  At two of the vertices, say $v_1$ and $v_2$, three of the emanating  edges lie in $D$ and one does not. Denote the ones that do not lie in $D$ by $e_1$ and $e_2$, respectively. At all other  $v_i$ we have two edges that lie in $D$ and two that do not. Denote the latter by $e_i, e_i'$. We have a total of $2(k-2)+2=2k-2$ edges not in $D$. The problem we have now is combinatorial: 

 \textit{Given the disk $D'=\overline{S^2-D}$ and the points $v_1,...,v_k\in\partial D'$, $k\geq 4$, is it possible to subdivide $D'$ by  $2k-2$ edges in  a way that  exactly one edge emanates from  both $v_1$ and $v_2$ and exactly two edges emanate from $v_3,...,v_k$  in such a way that no pair of edges intersect and every face on the subdivision of $D'$ is at least 3-sided (here we also consider sides coming from the boundary)?}

A simple argument will show that the answer to this question is negative. 
Orient the boundary of $D'$  counterclockwise. Starting at $v_1$, draw the edge $e_1$ emanating from it. The other endpoint of $e_1$ is some vertex $v_{i_1}$. Consider the vertices contained in the segment  $[v_1,v_{i_1}]\subset\partial D'$ in the given orientation. If there are no vertices at all, then we must have a  2-sided face, which is not possible. Therefore, by relabeling, we may assume $v_2$ is the the first vertex between $v_1$ and $v_{i_1}$. Observe that the edges emanating from $v_2$ are trapped between the edge $e_1$ and $\partial D'$. Draw an edge $e_2$ emanating from $v_2$ with the second endpoint $v_{i_2}$. It must be that $v_{i_2}$ also lies in $[v_1,v_{i_1}]$, or else we find a pair of intersecting edges. As above, there must be a vertex in the segment $[v_2,v_{i2}]$.  By repeating the above  argument eventually we find a 2-sided face,  which is not possible.   Therefore it must be that $V_1> B_1+W_1-2$.
\end{proof}

\section{ CO-FINALNESS}

In this section we provide a way of choosing the black or white faces on the
polyhedra $P_j$ along which it is reflected in such a way that the resulting family $\{M_j\}$ of manifolds is
cofinal. The main result of this section, Theorem \ref{G_ij}, appears as part of the proof of Theorem 2.2 of \cite{Ag}. We include a proof here for completeness.
To better describe this construction we  need to change notation slightly by
adding another index. 

Start with $P_1$ and relabel it $P_{11}$. Reflect
along a black face $f_{11}$ obtaining $P_{12}$. Let $\phi_{f_{11}}$ represent such reflection.  	Observe that if $f$ is adjacent to $f_{11}$, then $f\cup\phi_{f_{11}}(f)$ corresponds to  a single face on $P_{12}$. We call $f$ and $\phi_{f_{11}}(f)$ \textit{subfaces} of $f\cup\phi_{f_{11}}(f)$. Next reflect $P_{12}$ along a
white face $f_{12}$, which is also a face  of $P_{11}$ or contais a face of $P_{11}$ as a subface, obtaining $P_{13}$. We construct a
subcollection $P_{11},...,P_{1k_1}$ of  polyhedra such that 
\begin{itemize}
\item[(i)] If $P_{1j}$ is obtained from $P_{1(j-1)}$ by reflection along a white
(black) face  then  $P_{1(j+1)}$ is obtained from $P_{1j}$ by reflection along a
black (white) face.
\item[(ii)] Whenever possible, the face $f_{1j}$ must be a face of $P_{11}$ or contain a  face of $P_{11}$ as a subface.
\item[(iii)]  No faces of $P_{11}$ are subfaces of $P_{1k_1}$.   
\end{itemize}

\noindent
 Now set $P_{1k_1}:=P_{21}$.  
\noindent
Suppose $P_{n1}$ has been constructed.
Construct the subcollection  of polyhedra $P_{n1},...,P_{nk_n}$ such that 

\begin{itemize}
\item[(i)] The reflections were performed in a alternating fashion with respect to the color of the faces;
\item[(ii)] Whenever possible, the face $f_{nj}$ must be a face of $P_{n1}$ or contain a  face of $P_{n1}$ as a subface.
\item[(iii)]  No faces of $P_{n1}$ are subfaces of $P_{nk_n}$. 
\end{itemize}

Now set $P_{nk_n}:=P_{(n+1)1}$. Inductively we obtain a collection of polyhedra
$$P_{11},P_{12},...,P_{1k_1}:=P_{21},...,P_{2k_2}:=P_{31},...,P_{nk_n}:=P_{(n+1)1},...$$
satisfying (i), (ii) and (iii) above.

Let $G_{ij}$ be the reflection group of  $P_{ij}$ and let $M_{ij}$ be the cover of $M_{11}$ whose fundamental group is $\pi_1(M_{ij})=\pi_1(M_{11})\cap G_{ij}$. Co-finalness of the family $\{M_{ij}\longrightarrow M_{11}\}$ is an immediate consequence of

\begin{theorem}\label{G_ij}
Let $G_{ij}$ be as above. Then $\displaystyle\cap_{ij}G_{ij}=\{1\}$.
\end{theorem}

In order to prove this theorem we consider the base point for the
fundamental group of each $P_{ij}$ (viewed as orbifolds with their faces mirrored) to be the barycenter $x_0$ of $P_{11}$.  

\begin{proof}[Proof of Theorem]

Set $R_{ij}=\displaystyle\inf_{\gamma}\{\ell(\gamma)\}$, where $\gamma$ is an arc
with endpoints in faces (possibly edges) of $P_{ij}$ going through $x_0$. Note
that, by construction,  $\displaystyle\lim_{i\rightarrow\infty }R_{ij}=\infty$. 
For a non-trivial element  $g\in G_{11}$ set
$R_g=\displaystyle\inf_{[\alpha]=g}\{\ell(\alpha)\}$, where $\alpha$ is a loop in $P_{11}$ based at $x_0$ and $[\alpha]$ represents its homotopy class. 
Let $\alpha_g$ be a loop in $P_{11}$ based at $x_0$ such that $[\alpha_g]=g$ and $\ell(\alpha_g)\leq R_g+1$.

We claim that for sufficiently large $i$ one cannot have
$g\in G_{ij}$. In fact, if $\alpha_{ij}$ is any loop in $P_{ij}$ based at $x_0$, then  this loop bounces off faces of $P_{ij}$, yielding  an arc  $\gamma_{ij}$
throught $x_0$. Therefore $\ell(\alpha_{ij})\geq\ell(\gamma_{ij})\geq R_{ij}$. Since covering maps preserve length of curves, this  implies that if $i$ is large enough no such $\alpha_{ij}$ maps to $\alpha_g$. Thus it is not possible to find a loop representative for $g$ in $P_{ij}$.  
\end{proof}

\section{FINAL REMARKS}

\begin{question}\label{cofinal} Is it possible,  in our setting, to obtain a co-final tower of regular covers $\{M_j\longrightarrow M_1\}$ with positive rank gradient?
\end{question}

  A positive answer to this  would be very relevant, as  it  implies that Question \ref{cost} has a negative answer.  However, the tower constructed in Theorem \ref{E1s1} cannot consist of normal subgroups.  To see this we argue as follows: using the main theorem in   \cite{Ma} we can find  a sequence $\{\gamma_j\}$ of hyperbolic elements, $\gamma_j\in G_j$, whose translation lengths are bounded above by $2.634$. Since there exist at most finitely many conjugacy classes of hyperbolic elements of bounded translation length in $G_1$, it must be that an infinite subsequence $\{\gamma_{j_k}\}$  lie in the same conjugacy class in $G_1$.  Let $\gamma$ be a representative of this  class and $g_{j_k}\in G_1$ be such that $\gamma_{j_k}=g_{j_k}\gamma g_{j_k}^{-1}$. If the tower $\{G_j\}$ consists of  normal subgroups, then $\gamma\in G_{j_k}$, contradicting the fact that $\{G_{j_k}\}$ is co-final. 
  
Question \ref{cofinal} is relevant  also because of the following result  (see \cite{AN}):

\begin{thm}[Ab\'ert-Nikolov] Either the Rank vs. Heegaard genus conjecture (see below) is false or  Question \ref{cost}  has a negative solution.
\end{thm}

If an orientable 3-manifold M is closed, a Heegaard splitting of  $M$ consists of two handlebodies $H_1$ and $H_2$ with their boundaries identified by some orientation preserving homeomorphism.  Recall that the genus of, say, $\partial H_1$ gives an upper bound on the rank of $\pi_1(M)$. If $M$ is not closed, these decompositions are given in terms of compression bodies, again denoted by $H_1$ and $H_2$.  In order to obtain useful bounds on the rank of $\pi_1(M)$ we restrict ourselves to those decompositions in which $H_1$, for instance, is a  handlebody. Note that if this is the case, then the  genus of $\partial H_1$ is again an upper bound for the rank of $\pi_1(M)$. Recall that the \textit{Heegaard genus} of $M$ is the minimal genus of a Heegaard surface. A long standing question in 3-dimensional topology is:

\begin{conjecture} The rank of an orientable finite volume hyperbolic 3-manifold equals its Heegaard genus.
\end{conjecture}

Another concept due to Lackenby is that of \textit{Heegaard gradient} (\cite{La2}). 
Given a orientable 3-manifold $M$and a family $\{M_j\}$ of finite sheeted covers, we define the \text{Heegaard grandient of $\{M_j\longrightarrow M\}$} by 
$$\text{Hgr}(M,\{M_j\})=\lim_{j\rightarrow\infty}\frac{-\chi(S_j)}{d_j}$$   
where $d_j$ is the degree of the cover $M_j\longrightarrow M$ and $S_j$ is a minimal genus Heegaard surface for $M_j$. 

Note that if $\text{rgr}(M,\{M_j\})>0$, then $\text{Hgr}(M,\{M_j\})>0$.
An important conjecture that would follow from the  \textquotedblleft rank versus Heegaard genus\textquotedblright conjecture  is 

\begin{conjecture} Let $M$ be a finite volume hyperbolic 3-manifold and $\{M_i\longrightarrow M\}$ a family of finite sheeted covers. Then $\operatorname{rgr}(M,\{M_i\})>0$ if and only if $\operatorname{Hgr}(M,\{M_i\})>0$
\end{conjecture}

Our results provide examples for which this is true. 
In (\cite{La2}) Lackenby showed that if $\pi_1(M)$ is an arithmetic lattice in PSL(2,$\mathbb{C}$),  then M has a co-final family of covers (namely, those arising from congruence subgroups) with positive  Heegaard gradient.  In \cite{LLR} Long, Lubotzky and Reid generalize this result by proving that every finite volume hyperbolic 3-manifold has a co-final family of finite sheeted regular covers for which the Heegaard gradient is positive. These results were also motivation for this note. 

A natural question  that arises from our results is to what other categories of finite volume   hyperbolic 3-manifolds they hold. For instance: 

\begin{question} Is it true that given a right-angled poyhedron $P_1$ (not necessarily ideal) and  a manifold $M_1$ such that $\pi_1(M_1)$ has finite index in the reflection group of $P_1$, then there exists a co-final tower $\{M_j\longrightarrow M_1\}$ of finite sheeted covers with  positive rank gradient?
\end{question}
In our setting the ideal vertices played an important role as they were used to find lower bounds on the rank of the fundamental groups.  If the polyhedron $P_1$ has  vertices which are not ideal  then we need to find another way of estimating the rank of the associated manifolds. Ian Agol has suggested a way for  doing this. We are currently working on appropriate bounds for the rank in this case and will include it in a future work. 


It is also easy to give examples of families $\{M_j\longrightarrow M_1\}$ with arbitrarily large rank gradient. Using the methods above it suffices to provide examples of polyhedra $P_1$ for which the difference $V_1-(B_1+W_1)$ is arbitrarily large. Below we illustrate some cases in which this happens: consider the right-angled ideal polyhedron $P_0$ pictured below, viewed as lying in $S^2$. 
 
\begin{figure}[h]
\includegraphics[scale=.3]{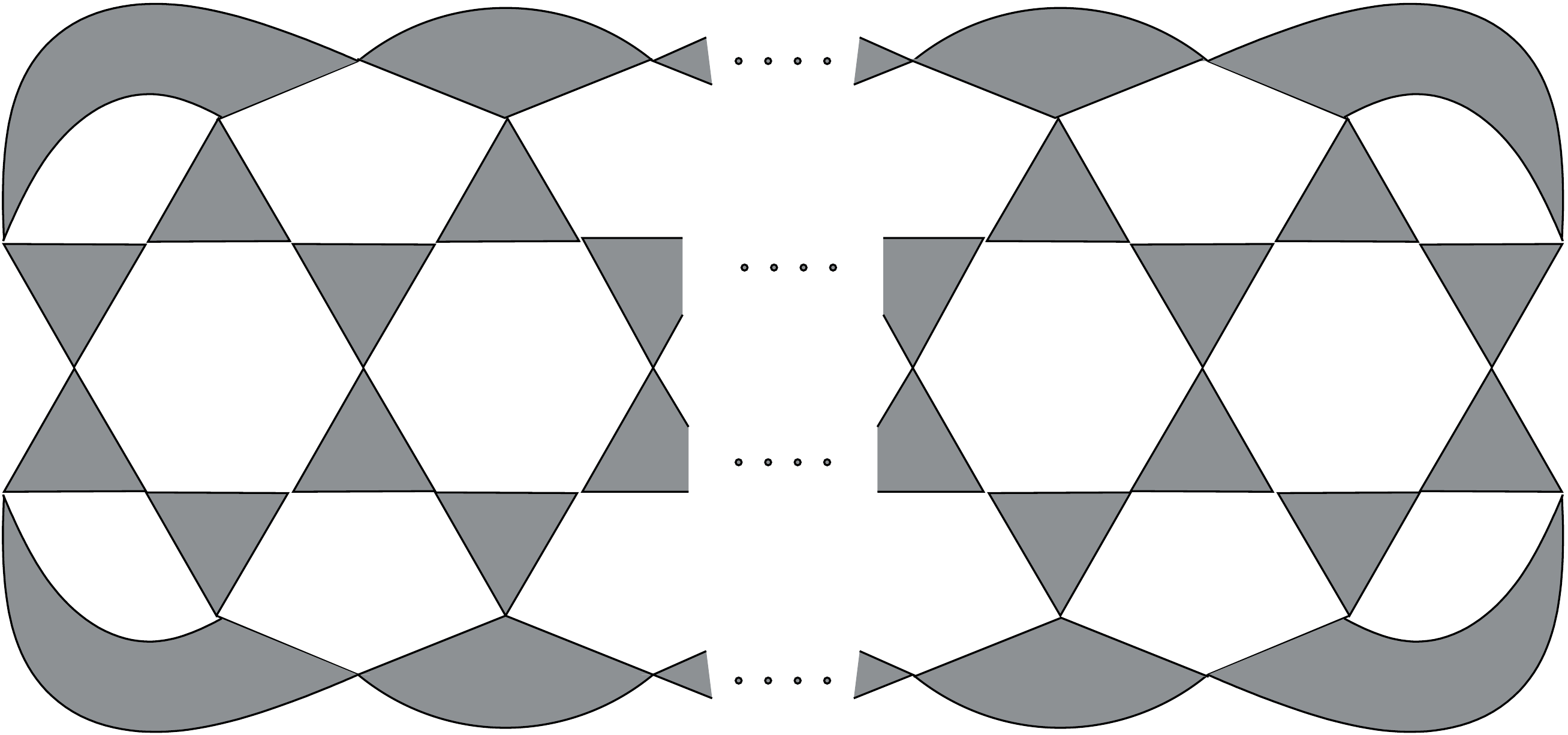}
\caption{Polyhedron $P_0$}
\end{figure}

Note that, by Andreev's theorem, this polyhedron can be realized as a totally geodesic right-angled ideal polyhedron in $\mathbb{H}^3$.  
Reflecting $P_0$ along the white face containing the point at infinity of $S^2$ will give us a polyhedron $P_1$. Since $P_1$ is obtained from  two copies of  $P_0$ by gluing together the white faces containing the point at infinity, we have a maximum of $6$ ideal vertices per white face of $P_1$ and a maximum of $4$ per black faces. Obviously  this construction can be made so that $P_1$ has arbitrarily many ideal vetices. Thus, given any $C>0$ we may find $P_1$ such that for the family $\{M_j\longrightarrow M_1\}$ as above 
$$\lim_{j\rightarrow\infty}\dfrac{\text{rk}(\pi_1(M_j))-1}{[\pi_1(M_1):\pi_1(M_j)]}\geq\lim_{j\rightarrow\infty}\dfrac{
2^{j-1}(V_1-(B_1+W_1))-1}{2^{j-1}}>C$$

\vspace{.4cm}
\noindent
\address{\textsc{Department of Mathematics,\\
The University of Texas at Austin}}\\
\email{\textit{E-mail:}\texttt{ dgirao@math.utexas.edu}}

\end{document}